\documentclass{article}
\usepackage{spconf}
\usepackage{cite}
\usepackage{amsmath,amssymb,amsfonts}
\usepackage{graphicx}
\usepackage{textcomp}
\usepackage{xcolor}
\usepackage{theorem}
\usepackage{charter}
\usepackage{mathrsfs} 
\usepackage{euscript}
\usepackage{enumitem}
\usepackage{caption}
\usepackage{multirow}
\usepackage{array}
\usepackage{booktabs} 
\usepackage[usenames,dvipsnames]{pstricks}
\definecolor{dred}{HTML}{D90404}

\usepackage{url}
\newcommand{\email}[1]{\href{mailto:#1}{\nolinkurl{#1}}}
\renewcommand{\leq}{\ensuremath{\leqslant}}
\renewcommand{\geq}{\ensuremath{\geqslant}}

\newcommand{\minimize}[2]{\ensuremath{\underset{\substack{{#1}}}%
{\text{\rm minimize}}\;\;#2 }}

\newcommand{\Scal}[2]{\bigg\langle{#1}\;\bigg|\:{#2}\bigg\rangle} 
\newcommand{\scal}[2]{{\left\langle{{#1}\mid{#2}}\right\rangle}}

\newcommand{\HH}{\ensuremath{{\mathcal H}}}
\newcommand{\GG}{\ensuremath{{\mathcal G}}}

\newcommand{\emp}{\ensuremath{{\varnothing}}}

\newcommand{\RR}{\ensuremath{\mathbb{R}}}

\newcommand{\NN}{\ensuremath{\mathbb N}}

\newcommand{\KK}{\ensuremath{\mathbb K}}

\newcommand{\thetab}{\ensuremath{\boldsymbol{\theta}}}  
\newcommand{\Thetab}{\ensuremath{\boldsymbol{\Theta}}}

\newcommand{\zeroun}{\ensuremath{\left]0,1\right[}}

\newcommand{\Argmind}[2]{\ensuremath{\underset{\substack{{#1}}}%
{\mathrm{Argmin}}\;\;#2 }}
\newcommand{\proj}[1]{\ensuremath{\text{\rm proj}_{#1}\,}}

\newtheorem{theorem}{Theorem}

\newtheorem{proposition}[theorem]{Proposition}
\theoremstyle{plain}{\theorembodyfont{\rmfamily}%
}
\theoremstyle{plain}{\theorembodyfont{\rmfamily}%
}
\theoremstyle{plain}{\theorembodyfont{\rmfamily}%
\newtheorem{algorithm}[theorem]{Algorithm}}
\theoremstyle{plain}{\theorembodyfont{\rmfamily}%
}
\theoremstyle{plain}{\theorembodyfont{\rmfamily}%
}
\theoremstyle{plain}{\theorembodyfont{\rmfamily}%
}
\theoremstyle{plain}{\theorembodyfont{\rmfamily}%
}
\theoremstyle{plain}{\theorembodyfont{\rmfamily}%
}

\title{A Variational Inequality Model for Learning Neural Networks}

\name{Patrick L. Combettes$^1$ \qquad 
Jean-Christophe Pesquet$^2$ \qquad 
Audrey Repetti$^3$ %
\thanks{The work of P. L. Combettes was supported by the 
National Science Foundation under grant CCF-2211123.
The work of J.-C. Pesquet was supported by
the ANR Chair in AI BRIDGEABLE. 
The work of A. Repetti was supported by an RSE Personal Research
Fellowship from the Royal Society of Edinburgh, by the CVN,
INRIA/OPIS and CentraleSup\'elec.}}
\address{\small$^1$ \textit{Department of Mathematics},
\textit{North Carolina State University},
Raleigh, USA\\
\small$^2$ \textit{Centre de Vision Num\'erique},
\textit{CentraleSup\'elec},
\textit{Inria, Unversit\'e Paris-Saclay},
Gif sur Yvette, France\\
\small$^3$ \textit{Maxwell Inst. for Mathematical 
Sciences, Inst. of Sensors, Signals and Systems},
\textit{Heriot-Watt University},
Edinburgh, UK
}

\begin{document}
\maketitle

\begin{abstract}
Neural networks have become ubiquitous tools for solving signal and
image processing problems, and they often outperform standard
approaches. Nevertheless, training neural networks is a
challenging task in many applications. The prevalent training
procedure consists of minimizing highly non-convex objectives
based on data sets of huge dimension. In this context, current
methodologies are not guaranteed to produce global solutions. We
present an alternative approach which foregoes the optimization
framework and adopts a variational inequality formalism. The
associated algorithm guarantees convergence of the iterates to a
true solution of the variational inequality and it possesses an
efficient block-iterative structure. A numerical application is
presented.
\end{abstract}

\begin{keywords}
Block-iterative algorithm, MRI, neural networks, 
transfer learning, variational inequality. 
\end{keywords}

\section{Introduction}
\label{sec:1}

Deep learning techniques have become very successful in solving a
great variety of tasks in data science; see for instance
\cite{Adle18,ahmad2020plug,Hayk98,Heka12,Lecu15,Zhan19,Zhan21}.
Deep neural networks rely on highly parametrized nonlinear systems.
Standard methods for learning the vector of parameters $\thetab$ of
a neural network ${T}_{\thetab}$ are mainly based on stochastic
algorithms such as the stochastic gradient descent (SGD) or Adam
methods \cite{Kingma2015,Zhan19}, and they are implemented in
toolboxes such as PyTorch or TensorFlow. In this context, the
standard approach to learn the parameter vector $\thetab$ is to
minimize a training loss. Specifically, given a finite training
data set consisting of ground truth input/output pairs
$(x_k,y_k)_{1\leq k\leq K}$, a discrepancy measure is computed
between the ground truth outputs $(y_k)_{1\leq k\leq K}$ and the
outputs $(T_{\thetab}x_k)_{1\leq k\leq K}$ of the neural network
driven by inputs $(x_k)_{1\leq k\leq K}$. Thus, if we denote by
$\Thetab$ the parameter space, the objective of these methods is to 
\begin{equation}
\label{e:1}
\minimize{\thetab\in\Thetab}\sum_{k=1}^K\ell\big(
T_{\thetab}x_k,y_k\big).
\end{equation}
One of the main weaknesses of such an approach is that it typically
leads to a nonconvex optimization problem, for which existing
algorithms offer no guarantee of optimality for the delivered output
parameters. In other words, the solution methods do
not provide true solutions but only local ones that may be hard to
interpret in terms of the original objectives in~\eqref{e:1}. 

The contribution of this work is to introduce an alternative
training approach which is not based on an optimization approach
but, rather, seeks the parameter vector $\thetab$ as the solution
of equilibrium problems defined by variational inequalities.
Nonlinear analysis tools for neural network modeling have been
employed in 
\cite{Ree18,Coh19,Bai19,Svva20,Smds20,Hasa20,Wins20,Xu20,Pesq21}.
Here, we show that training a layer of a feedforward neural
network can be modeled as a variational inequality problem and
solved efficiently via iterative techniques such as the
deterministic block-iterative forward-backward 
algorithm of \cite{Anon21}. This
algorithm displays two attractive features. First, it guarantees
convergence of the iterates to a true equilibrium, and not to a
local solution as in the minimization setting. Second, it lends
itself to an implementation based on a batch strategy, which is
indispensable to deal with large data sets. The strategy of
foregoing standard optimization in favor of more general forms of
equilibria in the form of variational inequalities was first
adopted in \cite{Siim22} in a quite different context, namely
signal recovery in the presence of nonlinear observations. 

The paper is organized as follows. Section~\ref{sec:FP}
describes our new training method, the design of a mini-batch
algorithm to solve the associated variational inequality, and
convergence properties of this algorithm. In Section~\ref{sec:TL},
we apply the proposed approach to a transfer learning problem in
which the last layer of neural network is optimized to denoise
magnetic resonance (MR) images. Some conclusions are drawn in
Section~\ref{sec:conclu}.

\section{Proposed variational inequality model}
\label{sec:FP}

\subsection{Variational inequality model for a single layer}
\label{sec:mo}

We first consider a single layer, modeled by an
operator $T_\theta$ acting between an input Euclidean space $\HH$
and output Euclidean space $\GG$, and parametrized by a vector
$\theta$ which is constrained to lie in a closed convex subset 
$C$ of a Euclidean space $\Theta$. More specifically,
\begin{equation} 
\label{e:NNlay}
T_\theta\colon\HH\to\GG\colon x\mapsto R(Wx+b),
\end{equation}
where $W\colon\HH\to\GG$ is a linear weight operator, 
$b\in\GG$ a bias vector, and
$R\colon\GG\to\GG$ a known activation operator. The objective is to
learn $W$ and $b$ from a training data set 
$(x_k,y_k)_{1\leq k\leq K}\in(\HH\times\GG)^K$. Our model
assumes that the parametrization $\theta\mapsto(W,b)$ is linear.
Further, we set
\begin{equation}
\label{e:5}
(\forall k\in\{1,\ldots,K\})\quad
L_k\colon \Theta\to\GG\colon\theta \mapsto  Wx_k+b.
\end{equation}
Thus, the ideal problem is to
\begin{multline}    
\label{e:2-}
\text{find}\;\;\theta\in{C}\:\;\text{such that}\\
(\forall k\in\{1,\ldots,K\})\:\;T_\theta x_k=y_k,
\end{multline}
that is, to
\begin{multline}    
\label{e:2}
\text{find}\;\;\theta\in{C}\;\:\text{such that}\\
(\forall k\in\{1,\ldots,K\})\;\:R(L_k\theta)=y_k.
\end{multline}
In practice, this ideal formulation has no solution and one must
introduce a meaningful relaxation of it. This is usually done via
optimization formulations such as \eqref{e:1}, which leads to the
pitfalls discussed in Section~\ref{sec:1}.

The approach we propose to construct a relaxation of \eqref{e:2}
starts with the observation made in \cite{Svva20} that most
activation operators are firmly nonexpansive in the sense that, for
every $z_1\in\GG$ and every $z_2\in\GG$, 
$\scal{z_1-z_2}{Rz_1-Rz_2}\geq\|Rz_1-Rz_2\|^2$.
Using this property, we can show that \eqref{e:2} can be relaxed
into the variational inequality problem 
\begin{multline}    
\label{e:3}
\text{Find}\;\;\theta\in C\;\text{such that}\\
(\forall\vartheta\in C)\;
\Scal{\vartheta-\theta}{\sum_{k=1}^K\omega_k L_k^*
\big(R(L_k\theta)-y_k\big)}\geq 0
\end{multline}
where, for every $k\in\{1,\ldots,K\}$, $L_k^*\colon\GG\to\Theta$ is
the adjoint of $L_k$ and $\omega_k\in\zeroun$, and
$\sum_{l=1}^K\omega_l=1$. This relaxation
is exact in the sense that, if \eqref{e:2} has solutions,
they are the same as those of \eqref{e:3} \cite{2023}. We
assume that \eqref{e:3} has solutions, which is true under mild
conditions \cite{2023}.

\subsection{Block-iterative forward-backward splitting}
\label{Ssec:FP:bFB}

We solve the variational inequality problem \eqref{e:3} 
by adapting a block-iterative forward-backward algorithm proposed
in \cite{Anon21}. This algorithm splits the computations associated
with the different linear operators $(L_k)_{1\leq k\leq K}$ using a
block-iterative approach. At iteration $n\in\NN$, a subset
$\KK_n$ of $\{1,\ldots,K\}$ is selected and, for every $k\in\KK_n$,
a forward step in the direction of the vector
$L_k^*(R\big(L_k\theta_n\big)-y_k)$ is performed. The forward steps
are then averaged and projected onto the constraint set $C$.

\begin{algorithm}
\label{algo:1}
Take $\gamma\in\left]0,2/\max_{1\leq k\leq K}\|L_k\|^2\right[$, 
$\theta_0\in\Theta$,
and $(\vartheta_{k,0})_{1\leq k\leq K}\in\Theta^K$.
Iterate
\begin{equation}
\label{algo:bFB-trainlast}
\begin{array}{l}
\text{for}\;n=0,1,\ldots\\
\left\lfloor
\begin{array}{l}
\displaystyle\text{select}\;\emp\neq\KK_n\subset\{1,\ldots,K\}\\
\displaystyle\text{for every}\;k\in\KK_n\\
\left\lfloor
\begin{array}{l}
\displaystyle {\vartheta}_{k,n+1}=\theta_n-\gamma L_k^*
\big(R(L_k\theta_n)-y_k\big)
\end{array}
\right.\\
\displaystyle\text{for every}\;k\in\{1,\ldots,K\}
\smallsetminus\KK_n\\
\left\lfloor
\begin{array}{l}
\displaystyle {\vartheta}_{k,n+1}={\vartheta}_{k,n}
\end{array}
\right.\\
\displaystyle\theta_{n+1}=\proj{C} 
\left(\sum_{k=1}^K\omega_k\vartheta_{k,n+1}\right).
\end{array}
\right.
\end{array}
\end{equation}
\end{algorithm}

We then derive the following result from \cite{Anon21}.

\begin{proposition}
\label{p:1}
Suppose that, for some $P\in\NN$, every index $k\in\{1,\ldots,K\}$
is selected at least once within any $P$ consecutive iterations,
i.e., $(\forall n\in\NN)$ 
$\bigcup_{k=0}^{P-1}\KK_{n+k}=\{1,\ldots,K\}$.
Then the sequence $(\theta_n)_{n\in\NN}$ generated by
Algorithm~\ref{algo:1} converges to a solution to \eqref{e:3}.
\end{proposition}

\subsection{Proposed deterministic batch algorithm}
\label{sec:a}

In a neural network context, batch approaches are necessary for
training purposes. Towards this goal, we modify 
Algorithm~\ref{algo:1} into a batch-based deterministic
forward-backward scheme to solve \eqref{e:3}.

Let us form a partition
$(\KK_j)_{1\leq j\leq J}$ of $\{1,\ldots,K\}$ and assume that, at
each iteration $n\in\NN$, only one batch index
$j_n\in\{1,\ldots,J\}$ is selected. Then, to avoid keeping in
memory all values of $(\vartheta_{k,n})_{1\leq k\leq K}$, 
Algorithm~\ref{algo:1} is rewritten below (Algorithm~\ref{algo:2})
so that only $J$ variables are kept in memory. Given
$j_n\in\{1,\ldots,J\}$, $\overline{\theta}_{j_n,n}\in\Theta$
denotes the stored variable associated with subset $\KK_{j_n}$. 

\begin{algorithm}
\label{algo:2}
Take $\gamma\in\left]0,2/\max_{1\leq k\leq K}\|L_k\|^2\right[$, 
$\theta_0\in\Theta$, and
$(\overline{\theta}_{j,0})_{1\leq j\leq J}\in\Theta^J$.
Iterate
\begin{equation}
\begin{array}{l}
\text{for}\;n=0,1,\ldots\\[0.1cm]
\left\lfloor
\begin{array}{l}
\displaystyle\text{select}\;j_n\in\{1,\ldots,J\}\\
\displaystyle\overline{\theta}_{j_n, n+1}
=\sum_{k\in\KK_{j_n}}\!\omega_k\Big(\theta_n -
\gamma L_k^*\big(R(L_k\theta_n)-y_k\big)\Big)\\
\displaystyle\text{for}\;j\in\{1,\ldots, J\} 
\smallsetminus\{j_n\}\\[0.1cm]
\left\lfloor
\begin{array}{l}
\displaystyle\overline{\theta}_{j, n+1} 
=\overline{\theta}_{j,n} 
\end{array}
\right.\\[0.3cm]
\displaystyle\overline{\theta}_{n+1} 
=\overline{\theta}_n 
-\overline{\theta}_{j_n, n+1} 
+\overline{\theta}_{j_n, n}\\[0.2cm]
\displaystyle\theta_{n+1}=\proj{C} 
\overline{\theta}_{n+1}.
\end{array}
\right.
\end{array}
\end{equation}
\end{algorithm}

The sequence $(\theta_n)_{n\in\NN}$ in \eqref{algo:2} converges to 
a solution to \eqref{e:3} provided that there exists $P\in\NN$ such
that $(\forall n\in\NN)$ 
$\bigcup_{k=0}^{P-1}\{j_{n+k}\}=\{1,\ldots,J\}$ \cite{2023}.

\subsection{The case of general feedforward neural networks}
\label{sec:feedNN}

Let $\HH_0,\ldots,\HH_M$ be Euclidean spaces. A feedforward neural 
network ${T}_{\thetab}\colon\HH_0\to\HH_m$ consists of a
composition of $M$ layers 
\begin{equation}
\label{e:NN-layers} 
T_{\thetab}=T_{M,\theta_M}\circ
\cdots\circ T_{1,\theta_1}
\end{equation} 
where the operators $(T_{m,\theta_m})_{1\leq m\leq M}$ are as in 
Section~\ref{sec:mo}:
$\theta_m\in\Theta_m$ is a vector linearly parametrizing a
weight operator $W_m\colon\HH_{m-1}\to\HH_m$ and a bias vector
$b_m\in\HH_m$, and $R_m\colon\HH_m\to\HH_m$ is a firmly
nonexpansive activation operator. For convenience, we gather 
the learnable parameters of the network in a vector 
$\thetab=(\theta_m)_{1\leq m\leq M}$. Given a training sequence 
$(x_k,y_k)_{1\leq k\leq K}\in(\HH_0\times\HH_M)^K$, the approach
proposed in Section~\ref{sec:FP} is used to train the last layer of
the neural network. For this layer the input sequence is defined 
by $(\forall k\in\{1,\ldots,K\})$ 
$\widetilde{x}_k=T_{M-1,\theta_{M-1}}\circ
\cdots \circ T_{1,\theta_1} x_k$. 
Two learning scenarios are of special interest:
\begin{itemize}
\item 
Greedy training \cite{Beli2020}. 
Layers are added one after the other to form a
deep neural network.
\item 
Transfer learning \cite{Do2005,Bane07,Bird20,Weis16}.
The goal is to retrain the last layer of a trained neural network
to allow it to be applied to a different data set or a different
task. 
\end{itemize}

\section{Transfer learning: Fine-tuning last layer of a denoising
neural network}
\label{sec:TL}

We apply the proposed variational inequality model to a transfer
learning problem for building a denoising neural network. Transfer
learning \cite{Weis16,Do2005,Bane07,Bird20} is often used in
practice to tailor a neural network that has been trained on a
particular data set, for a different type of data and improve its
performance \cite{Weis16}. 

\subsection{General setting for denoising neural networks}
\label{Ssec:TL:gen}

We consider a denoising neural network 
${T}_{\thetab}\colon\HH\to\HH$ with $M$ layers, defined as
in \eqref{e:NN-layers}.
${T}_{\thetab^*}$ has been pretrained as a denoiser, such
that
\begin{equation}
\thetab^*\in\Argmind{\thetab\in\Thetab}
\sum_{k=1}^{K'}\ell\big({T}_{\thetab}u_k,v_k\big),
\end{equation}
where each $(u_k,v_k)\in\HH^2$ is a pair of
noisy/ground truth images, and $\ell\colon\HH\times\HH\to\RR$ a
loss function. The objective is to retrain only the last layer of
${T}_{\thetab^*}$ in order to use it on a different type of
images. For instance, if the network has been trained on natural
images, it can be fine-tuned to denoise medical images obtained by
modalities such as MR or computed tomography.

\subsection{Simulation setting}

In our experiments, ${T}_{\thetab^*}$ is a DnCNN with $M=20$
layers, of the form of \eqref{e:NN-layers}, where
$\HH_0=\RR^{N\times N}$, $\HH_1=\cdots=\HH_{19}=\RR^{64\times
N\times N}$, and $\HH_{20}=\RR^{N \times N}$. The layers of the
networks are as follows. For the first layer, $W_1$ represents a
convolutional layer with $1$ input, $64$ outputs, and a kernel of
size $3\times 3$. For every $m\in\{2,\ldots,19\}$, $W_m$
represents a multi-input multi-output convolutional layer with $64$
inputs, $64$ outputs, and a kernel of size $3\times 3$. Finally,
$W_{20}$ represents a convolutional layer with $64$ inputs, $1$
output and a kernel of size $3\times 3$. We use LeakyReLU
activation functions with negative slope parameter $10^{-2}$. As
shown in \cite{Svva20}, this operator is firmly nonexpansive. In
addition, we take $b_1=\cdots=b_{20}=0$.

The network ${T}_{\thetab^*}$ is
trained on the $50,000$ ImageNet test data set converted to
grayscale images and normalized between $0$ and $1$. The ground
truth images $(v_k)_{1 \leq k \leq K'}$ correspond to patches of
size $50 \times 50$ selected randomly during training. For every
$k\in\{1,\ldots,K'\}$, the degraded images are obtained as
$u_k=v_k+\sigma b_k$, where $\sigma=0.02$ and $b_k\in\RR^{50\times
50}$ is a realization of a random standard white Gaussian noise.

We propose to fine-tune the network ${T}_{\thetab^*}$ to denoise MR
images. We thus focus on the training of the last layer 
$T_{20,\theta_{20}}$. We
choose $W_{20}$ to be a convolutional layer with 
$64$ inputs, 
$1$ output,
and kernels $w$ of size $7\times 7$. In addition, $R_{20}$
is a LeakyReLU activation function with negative slope
parameter $10^{-3}$. In this context, \eqref{e:5} assumes the form 
\begin{equation}
\label{e:Lj-conv}
L_k\colon\RR^{64\times 1\times 7\times 7}\to\RR^{N\times N} 
\colon w\mapsto\widetilde{x}_k * w,
\end{equation}
where $\widetilde{x}_k\in\RR^{64\times N\times N}$ is the
output of the 19th layer.

Three training strategies for $T_{20,\theta_{20}}$ are considered:
the standard SGD and the Adam algorithm for minimizing an $\ell^1$ 
loss, as well as the approach proposed in Section~\ref{sec:FP} with 
$C=\Theta_{20}=\RR^{64\times 1\times 7\times 7}$.

For the three methods, the
training set consists of the first $300$ images of the fastMRI
train data set, and we test the resulting networks on the next
$300$ images of the same data set. The dimension of the images in
this data set is $320\times 320$. We train the networks on patches,
by splitting the images into $16$ patches of size $80\times 80$.
The patches are randomly shuffled every time the algorithm has seen
all the patches of the data set. Moreover, we split the training
set into batches containing $10$ patches located at the same
position in $10$ images of the train set. Batches are normalized
between $0$ and $1$, and corrupted with an additive white Gaussian
noise with standard deviation $0.07$. One epoch is completed when
the algorithm has seen all the batches at the same location (i.e.,
$30$ batches generated as explained above). All algorithms are run
over $1000$ epochs.

The learning rate in SGD and Adam has been tuned manually to reach
the best performances. For Algorithm~\ref{algo:2}, the best
performances are achieved by choosing a large parameter, namely
$\gamma=1.9/\max_{1\leq k\leq K}\|L_k\|^2$. 

\begin{figure}
\centering
\includegraphics[width=0.95\columnwidth,%
trim={0.4cm 0.4cm 0.5cm 0.5cm},clip]%
{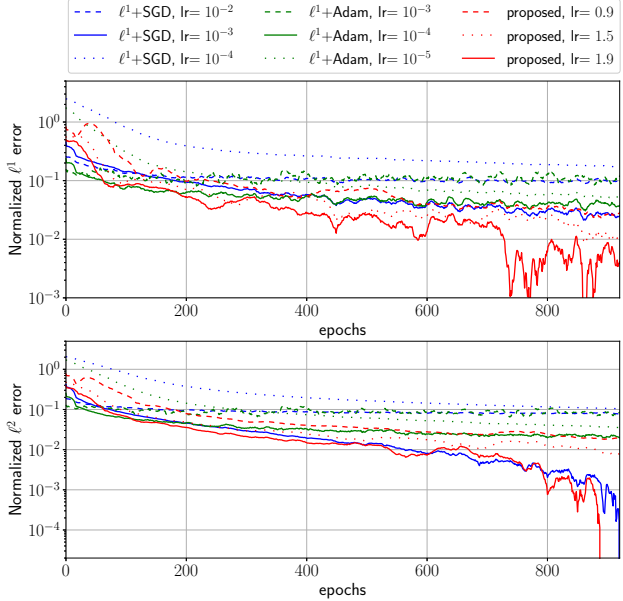}
\caption{\label{fig:cvgce}
\small
Convergence profiles showing the normalized averaged $\ell^1$ (top)
and $\ell^2$ (bottom) errors (in log scales) with respect to
epochs, for the $\ell^1$+SGD method (blue), the $\ell^1$+Adam
method (green), and the proposed approach (red). Continuous lines
show best step-size (i.e., learning rate) for each method. Dashed
and dotted lines show inaccurate choice of step-size.}
\end{figure}

\begin{table}[!]
\caption{
\label{tab:results}
\small
Average SSIM (and standard deviation) obtained 
for the first 300 images of the fastMRI training set,
and the next 300 images of the same set.}
\centering\small
\begin{tabular}{c|rr}
\hline
&\multicolumn{2}{c}{SSIM}\\
method & Training set & Test set \\
\hline
proposed &\textbf{0.6647 ($\pm$0.0721)}&\textbf{0.6630
($\pm$0.0597)}\\
$\ell^1$+SGD   &0.6641 ($\pm$0.0770) & 0.6627 ($\pm$0.0629)\\
$\ell^1$+Adam  &0.6598 ($\pm$0.0703) & 0.6239 ($\pm$0.0346)\\
\hline
\end{tabular}
\end{table}

\begin{figure}[t!]
\centering
\includegraphics[width=0.95\columnwidth,%
trim={2.6cm 3.2cm 1.1cm 1.5cm},clip]%
{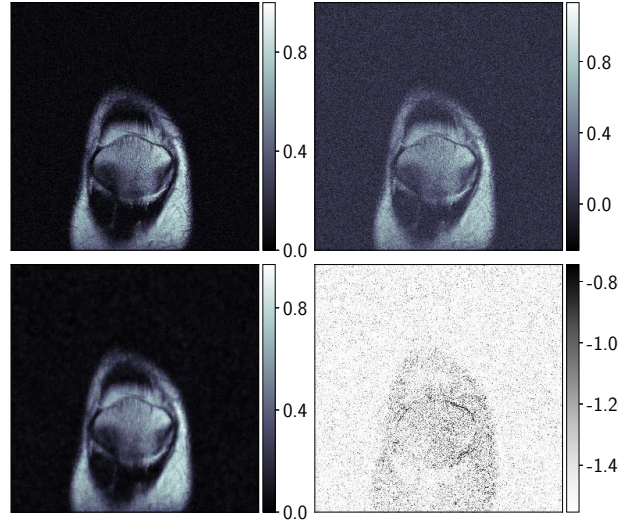}
    
\vspace{0.1cm}

\caption{\label{fig:res336}\small
Denoising results on the test set corresponding to slice 399
of the fastMRI training set. 
First row: ground truth (left) and corresponding noisy
observation (right) with PSNR $=23.09$ dB.
Second row: output of the DnCNN trained with the proposed
procedure (left), with PNSR $=29.31$ dB, and the
corresponding error map, in log-scale (right).
}
\end{figure}

\subsection{Simulation results}
Since the proposed method is not devised as a minimization 
method, we assess the behavior of the three learning procedure
during training by monitoring the $\ell^p$ errors 
$\sum_{k=1}^K\| T_{20,\theta_{20,e}}\widetilde{x}_k-y_k\|^p_p$
for $p\in\{1,2\}$ with respect to the epochs 
$e\in\{1,\ldots,1000\}$. 
We observe that, for any choice of the
step-size value $\gamma$ (even not determined optimally), our
method reaches a lower $\ell^1$ error more quickly than SGD and 
Adam, for any choice of the learning rate. 
For the $\ell^2$ error, any choice of
step-size will lead to faster convergence than Adam. For this
example, an accurate choice of learning rate for SGD leads to
a performance which is similar to that of the proposed approach.
However, choosing an inaccurate learning rate results in 
extremely slow convergence (to a local solution) or diverging
behavior for SGD, while our method converges to a true solution of
\eqref{e:3} for any choice of step-size as long as it satisfies the
conditions given in Algorithm~\ref{algo:2}.

The SSIM values for the 300 training images and the 300 test images
are shown in Table~\ref{tab:results}. We observe that
our approach yields slightly better results for both data
sets. One image slice of the test data set is displayed in
Fig.~\ref{fig:res336} to show the good visual quality of the
proposed transfer learning approach. 

\section{Conclusion}
\label{sec:conclu}

A new framework has been proposed to train neural network layers
based on a variational inequality model. The effectiveness of this
approach has been illustrated through simulations on a transfer
learning problem. In future work, we plan to explore further
algorithmic developments and consider various applications of the
proposed technique to other training problems.

\end{document}